\documentclass{amsart}

\usepackage{amsfonts}
\usepackage{cases}
\usepackage{mathrsfs}
\usepackage{bbm}
\usepackage{amssymb}
\usepackage{txfonts}
\usepackage{amscd}
\usepackage{amsfonts,latexsym,amsmath,amsthm,amsxtra,mathdots}
\usepackage[bookmarks,colorlinks]{hyperref}
\usepackage[all,cmtip]{xy}
\RequirePackage{amsmath} \RequirePackage{amssymb}
\usepackage{color}
\usepackage{colordvi}
\usepackage{multicol}
\usepackage{tikz}
\usepackage{hyperref}
\usepackage{cleveref}
\usepackage{graphicx}
\usepackage{amsmath}
\usepackage{amsmath,amscd}
\usepackage[normalem]{ulem}
\usetikzlibrary{matrix,arrows}
\usepackage{textcomp}
\usepackage{tikz-cd}
\usepackage{nicefrac}

\newenvironment{named}[1]
  {\def\namedthmname{#1}%
   \refstepcounter{namedthm}%
   \namedthm\def\@currentlabel{#1}}
  {\endnamedthm}
\makeatother

     \newcommand{\BD}{{\mathbb {D}}}

    \newcommand{\BQ}{{\mathbb {Q}}} \newcommand{\BR}{{\mathbb {R}}}
     \newcommand{\BT}{{\mathbb {T}}}

     \newcommand{\BZ}{{\mathbb {Z}}}
     \newcommand{\CB}{{\mathcal {B}}}

\def\-{^{-1}}

\newcommand{\delete}[1]{}

    \theoremstyle{plain}

\newtheorem{thm}{Theorem}[section]
\newtheorem{lem}[thm]{Lemma}
\newtheorem{prop}[thm]{Proposition}
\newtheorem{cor}[thm]{Corollary}
\newtheorem{rem}[thm]{Remark}

\newtheorem{question}[thm]{Question}

\newtheorem*{namedthm}{\namedthmname}
\newcounter{namedthm}

    \numberwithin{equation}{section}

\def\Proof{\noindent{\bf Proof}\quad}
\def\qed{\hfill$\square$\smallskip}

\begin{document}

\title{Riemannian foliation with exotic tori as leaves}


\author{F. Thomas Farrell}
\address{Department of Mathematical Sciences and Yau Mathematical Sciences
Center, Tsinghua University, 100084, Beijing, China}
\email{farrell@math.binghamton.edu}

\author{Xiaolei Wu}
\address{University of Bonn, Mathematical Institute, Endenicher Allee 60, 53115, Bonn, Germany}
\email{xwu@math.uni-bonn.de}


\date{December, 2018}


\begin{abstract}
We construct smooth fiber bundles such that the fibers are exotic tori and the total space has finite abelian fundamental group. This gives examples of a Riemannian foliation on a closed manifold whose leaves are exotic tori and whose total space has finite abelian fundamental group.
\end{abstract}

\maketitle

\section{introduction}

Recall that a (singular) Riemaninan foliation is  a (singular) foliation on a complete Riemannian manifold with the property that every geodesic that is perpendicular at one point to a leaf remains perpendicular to every leaf it meets (see \cite[Section 2]{GGR15} for the precise definitions and examples). And a $B$-foliation is just a singular Riemaninan foliation whose leaves are homeomorphic to some Bieberbach  manifolds. The following question was raised by Fernando Galaz-Garcia and Marco Radeschi in \cite[Introduction, p.604]{GGR15}

\begin{question} \label{RSF}
Does there exist a non-trivial example of a singular Riemanian foliation whose leaves are exotic tori? Moreover, can one find $B$-foliations by exotic tori  on simply connected manifolds and, in particular, on spheres?
\end{question}

In this note, we are not able to answer these questions completely, but we do find many non-trivial examples of    Riemanian foliations whose leaves are exotic tori. Our main theorem is as follows.

\begin{thm}
Let $\BT$ be an exotic  torus  of dimension $n\geq 12$ and $M$ a smooth simply connected closed $4$-manifold such that $H_2(M,\BZ)\cong \BZ^n$. Then there exists a smooth fiber bundle
$$ \BT \hookrightarrow  E \xrightarrow{p} M$$
such that the total space $E$ has finite abelian fundamental group.
\end{thm}

\begin{rem} \label{exo-fiber}
We do not know any smooth fiber bundle $F \hookrightarrow E\rightarrow M$ of closed manifolds such that $F$ is an exotic torus, $E$ and $M$ are smooth simply connected manifolds. Note that the universal cover of the total space $E$ in our theorem also admits a fiber bundle structure whose fiber now is a finite cover of the exotic torus $\BT$. However, we do not know whether the new fiber is still an exotic torus as we do not have much control over the covering.
\end{rem}

It is not hard to see that given any smooth fiber bundle of smooth closed manifolds, we can put  Riemannian metrics on its total and base spaces such that the projection map is a Riemannian submersion (see \cite[Chapter 2]{GW09}, or  \cite[Section 17.9, p.205]{Mi08}). Since any Riemannian submersion gives rise to a  Riemannian foliation \cite[Corollary 26.12]{Mi08}, we have the following corollary.

\begin{cor}
For any $m\geq 16$, there exists a closed Riemannian manifold $E$ of dimension $m$ such that $E$ admits a   Riemannian foliation whose leaves are exotic tori and $\pi_1(E)$ has finite abelian fundamental group. 
\end{cor}

\begin{rem}
Note that by \cite[Section 15A]{Wa70}, there is an exotic torus in every dimension $n>4$.
\end{rem}

All manifolds in this note are assumed to be closed and smooth.

\textbf{Acknowledgements.} The second author was partially supported by Prof. Wolfgang L\"uck's ERC Advanced Grant “KL2MG-interactions”
(no. 662400). He also wants to thank the first author for inviting him to visit the Yau Mathematical Sciences Center at Tsinghua University where part of this work was done. Special thanks go to Mauricio Bustamante for helping him fill out the details for the proof of Theorem \ref{htpg-top-m-diff}. We also would like to thank the anonymous referee for some helpful comments.

\section{Some  calculations in homotopy groups}
In this section, we prove some results on homotopy groups that we will need later in constructing the bundles. Given a smooth closed manifold $M$, let $\text{Diff}(M)$ denotes the diffeomorphism group of $M$, $\text{Top}(M)$ the homeomorphism group of $M$ and $\text{G}(M)$ the monoid of self homotopy equivalences of $M$. Respectively, let $\text{BTop}(M)$, $\text{BDiff}(M)$ and $\text{BG}(M)$ be the corresponding classifying space. Moreover, let $\text{Top}_0(M)$ (resp. $\text{Diff}_0(M)$) be the subgroup of $\text{Top}(M)$ (resp. $\text{Diff}(M)$) consisting of all those homeomorphisms (resp. diffeomorphisms) of $M$ which are homotopic to the identity map. 

The following two results can be found for example in \cite[Lecture 5]{Fa02}. For more general results, see \cite[Section 4]{ELPUW}.

\begin{thm}\label{htp-g-Top-Diff}
Let $\BT$ be a smooth manifold which is homeomorphic to the the $n$-dimensional torus with $n\geq 10$, then for $1 \leq i \leq \min \{\frac{n-7}{2},\frac{n-4}{3}\} $ we have 

$$\pi_i(\text{Top}(\BT)) \otimes \BQ \cong \begin{cases} \BQ^n & i=1; \\
                     0 &   i\geq 2.
       \end{cases}$$
$$\pi_i(\text{Diff}(\BT)) \otimes \BQ \cong \begin{cases} \BQ^n & i=1 ;\\
                   \oplus_{j=1}^{\infty} H_{i+1-4j}(\BT,\BQ) &   i\geq 2 \text{ and } n \text{ is odd}; \\
                     0 &   i\geq 2 \text{ and } n \text{ is even.}
       \end{cases}$$       

\end{thm}

\begin{lem} \label{htp-g-g}
Let $M$ be an aspherical manifold and $G(M)$ be the space of homotopy equivalences of $M$, then 
$$\pi_i(G(M))  \cong \begin{cases} Out(\pi_1(M)) & i=0 ;\\
                    \text{Center}(\pi_1(M)) &   i=1; \\
                     0 &   i\geq 2. 
       \end{cases}$$  
\end{lem}

The following theorem might be well-known to the experts, but we did not find it in the literature. 

\begin{thm}\label{htpg-top-m-diff}
Let $\BT$ be a smooth manifold which is homeomorphic to  the $n$-dimensional torus with $n\geq 12$. Then $\pi_i( \text{Top}_0(\BT) /\text{Diff}_0(\BT)) \otimes \BQ \cong \{0\}$ for $i=1,2$.
\end{thm}

\begin{named}{Remark 2.3a}
For a nilpotent group $G$, $G\otimes \BQ = 0$ means that $G$ is a torsion group.
\end{named}

\Proof All the ingredients of the proof can be found in \cite{BFJ}. We first need Morlet's comparison theorem  (see \cite{BL74} or \cite[Section 2]{BFJ}). Let  $\text{Top}(n) = \text{Top}(\BR^n)$ and $O(n)$ be the $n$-th orthogonal group. Given a smooth manifold $M$, the tangent bundle of $M$, regarded as an Euclidean vector bundle, has an associated (right) principal $O(n)$-bundle $PM\rightarrow M$. Note that $O(n)$ acts on the left on the coset space $\text{Top}(n)/O(n)$. We can form the balanced product
$$ \CB_n(M):PM\times_{O(n)} \text{Top}(n)/O(n)$$
which is a fiber bundle over $M$ with fiber $\text{Top}(n)/O(n)$.  The space of sections of $\CB_n(M)$ is denoted by $\Gamma (\CB_n(M))$. Now Morlet's comparison theorem says, as long as $dim(M)\neq 4$, there exists a map $\text{Top}_0(M)/\text{Diff}_0(M) \rightarrow \Gamma (\CB_n(M))$ which induces an injective correspondence of connected components and a weak homotopy equivalence on each component. Now we take $M$ to be the exotic torus $\BT$, we have for $i\geq 1$
$$ \pi_i(\text{Top}_0(\BT)/\text{Diff}_0(\BT)) \cong \pi_i(\Gamma(\CB_n(\BT)))$$
 As indicated at the end of \cite[Section 4]{BFJ}, the  fiber $\text{Top}(n)/O(n)$ is simply connected and $\pi_i (\text{Top}(n)/O(n)) \otimes \BQ= 0$ for $0\leq i\leq n+2$ as $n\geq 12$. Note that here we do not need to assume  $n$ is odd as $\pi_i(\BD^{n},\partial) \otimes \BQ = 0$  for $n$ even in the Igusa stable range \cite{FH78} (see also \cite[Theorem 4.1]{RW17}). Now we apply \cite[Lemma 3]{BFJ} by taking both  fibration  as $\CB_n(\BT)$, we have $\pi_i(\Gamma(\CB_n(\BT)))\otimes \BQ = 0$ for $1 \leq i \leq 2$. Hence we have $\pi_i(\text{Top}_0(\BT)/\text{Diff}_0(\BT)) \otimes \BQ\cong 0$ for $1\leq i \leq 2$.

\qed

\begin{cor}\label{top-diff-iso}
Let $M$ be a smooth manifold which is homeomorphic to the the $n$-dimensional torus with $n\geq 12$. Then the  forgetful map $f: \text{BDiff}(\BT)\rightarrow \text{BTop}(\BT)$ induces a rational isomorphism in $\pi_2$.
\end{cor}

\Proof 
By the long exact sequence of homotopy groups associated to the fibration 

$$\text{Top}(\BT) /\text{Diff}(\BT)\rightarrow \text{BDiff}(\BT)\rightarrow \text{BTop}(\BT)$$
we only need to show $\pi_i( \text{Top}(\BT) /\text{Diff}(\BT)) \otimes \BQ \cong \{0\}$ for $i =1$ and $2$. This is proved in Theorem \ref{htpg-top-m-diff} since $\pi_i( \text{Top}(\BT) /\text{Diff}(\BT)) \cong \pi_i( \text{Top}_0(\BT) /\text{Diff}_0(\BT))$ for $i >0$.

\qed

\begin{cor}\label{loop-in-diff}
Let $\BT$ be an exotic torus of dimension $n \geq 12$ and $I:\text{BDiff}(\BT)\rightarrow \text{BG}(\BT)$ be the forgetful map. Then we can find based maps $f_i\in \pi_2(\text{BDiff}(\BT))$ for $1\leq i\leq n$, such that $I_\ast(f_i)$ generates a finite index subgroup in $\pi_2(\text{BG}(\BT))$.
\end{cor}
\Proof Note first that $I$ is the composition of the following maps: $$\text{BDiff}(\BT) \xrightarrow{\bar{i}_1} \text{BTop}(\BT) \xrightarrow{\bar{i}_2} \text{BG}(\BT)$$

Now to calculate the homotopy groups of $\text{BTop}(\BT)$, we do not need the smooth structure of $\BT$. Let $T$ be the standard torus which $\BT$ is homeomorphic to. Note that  $\pi_2(\text{BG}(T)) \cong \pi_1(\text{G}(T)) \cong\BZ^n$, and $\pi_2(\text{BTop}(T)) \otimes\BQ \cong \pi_1(\text{Top}(T)) \otimes\BQ \cong \BQ^n$. Thinking of $T$ as the $n$-dimensional flat torus ($\BR^n /\BZ^n$), then we can realize the generators of $\pi_1(\text{G}(T))$ as rotation along each factor  $\BR/\BZ$ in  $\BR^n/\BZ^n$. In particular they lies in $\pi_1(\text{Top}(T))$. Let $L$ denote the corresponding subgroup they generate in $\pi_2(\text{BTop}(T)) \cong \pi_1(\text{Top}(T))$. $L$ is a free abelian group of rank $n$.

By Corollary \ref{top-diff-iso}, we can find maps $f_i\in \pi_2(\text{BDiff}(\BT))$ such that ${\bar{i}_{1\ast}} (f_i)$ lies in $L$ and  generates a finite index subgroup in $L$. Hence $I(f_i)$ generate a finite index subgroup in $\pi_2(\text{BG}(\BT))$.

\qed

The following proposition also plays an important role in our construction of the bundle.

\begin{prop}\label{modify-nh}
Let $S^3 \xrightarrow{g} \bigvee_{i=1}^{k} S_i^2 \xrightarrow{f} X$ be two maps between spaces, where $S^3$ is the $3$-sphere and $S_i^2$ are $2$-spheres.  Suppose $\pi_3(X)\otimes\BQ \cong \{0\}$ and let $d_{i,N}$ be the degree $N$ self-map of $S_i^2$ which fixes the common wedge point. Then there exists $N$ such that  $f\circ d_N\circ g$ is homotopic to the constant map, where $d_N$ is the unique map such that $d_N\mid_{S^2_i} = d_{i,N}$.
\end{prop}

\Proof Note that in \cite[Theorem A]{Hi55}, Hilton calculated the homotopy groups for wedges of spheres. In our case, we have
$$\pi_3( \bigvee_{i=1}^{k} S_i^2) \cong \oplus_{i=1}^k\pi_3(S_i^2) \oplus (\oplus_{j=1}^{\frac{k(k-1)}{2}} \pi_3(S^3))$$
where $\pi_3(S^3)$ embeds in $\pi_3( \bigvee_{i=1}^{k} S_i^2)$ by  composing it with $S^3 \rightarrow S_{i_1}^2\bigvee S_{i_2}^2 \hookrightarrow \bigvee_{i=1}^{k} S_i^2$, $1\leq i_1<i_2\leq k$. Here the map from $S^3 \rightarrow S_{i_1}^2\bigvee S_{i_2}^2 $ is just the attaching map used in defining the Whitehead product. Now our map $g$ considered as an element in $\pi_3( \bigvee_{i=1}^{k} S_i^2)$ can be written as $x +y$, where $x\in \oplus_{i=1}^k\pi_3(S^2)$ and $y\in \oplus_{j=1}^{\frac{k(k-1)}{2}} \pi_3(S^3)$. By a  result of Hopf (see for example \cite[Chapter XI. Exercise A, Problem 4 and 5]{Hu59}), we have $d_{N,i} \circ x = N^2 x$ . On the other hand, since the Whitehead product is bilinear on its factors, we have $d_{N,i}\circ y =   N^2y$. Therefore,   $ d_N \circ g = N^2 g$.

Now since $\pi_3(X)\otimes\BQ \cong \{0\}$, we also have $ f \circ g$  is of finite order. Let $N$ be the order of $f\circ g$,  we have $f\circ d_N \circ g = f \circ (N^2 g) = N^2 f\circ g=0 \in \pi_3(X)$.

\qed

\section{Proof of the main theorem}
In this section, we prove our main theorem.

Recall that $M$ is a simply connected $4$-manifold such that $H_2(M,\BZ) \cong \BZ^n$ with $n\geq 12$ and $\BT$ is an exotic torus of dimension $n$. Note that a smooth bundle over $M$ with fiber the $n$-dimensional exotic torus $\BT$ is determined by a map $F: M \rightarrow \text{BDiff}(\BT)$ up to homotopy. We also have the following forgetful maps 

$$\text{Diff}(\BT) \xrightarrow{i_1} \text{Top}(\BT) \xrightarrow{i_2} \text{G}(\BT)$$

This in turn induces maps 

$$\text{BDiff}(\BT) \xrightarrow{\bar{i}_1} \text{BTop}(\BT) \xrightarrow{\bar{i}_2} \text{BG}(\BT)$$

Note that the fundamental group of the total space $E$ is determined by the composition $ \tilde{F} = \bar{i}_2\circ \bar{i}_1 \circ F: M \rightarrow \text{BG}(\BT)$. By Corollary  \ref{loop-in-diff}, we can find maps $f_i\in \pi_2(\text{BDiff}(\BT))$ such that $i_{2\ast}'\circ (i_{1\ast}') (f_i)\in  \pi_2(\text{BG}(\BT))$ ( $1\leq i \leq n$) generate a finite index subgroup in $\pi_2(\text{BG}(\BT))$.

Now we can construct the  map $F$ from $M$ to $\text{BDiff}(\BT) $ as follows. Note that since our $M$ is a  simply connected  $4$-manifold, it has a cell decomposition (up to homotopy equivalence, see the proof of \cite[Theorem1.2.25]{GS99}) with only one 0-cell, $n$ 2-cells and one 4-cell, i.e. M is  homotopy equivalent to  a cell complex.

$$M ~ = ~\bigvee_{i=1}^{n} S^2~ \bigsqcup_{S^3} ~D^4$$

We denote the attaching map $S^3 = \partial D^4 \rightarrow \bigvee_{i=1}^{n} S^2$ by $g$. We first define the map $F$ on $\bigvee_{i=1}^{n} S^2$ using the maps  $f_i\in \pi_2(\text{BDiff}(\BT))$. We are left to extend the map $F$ to $D^4$. Since  $F\mid_{\partial D^4 =S^3}$ is already defined, the map is extendable if and only if $F\mid_{\partial D^4 =S^3}$ represents a trivial element in $\pi_3(\text{BDiff}(\BT))$. But $\pi_3(\text{BDiff}(\BT)) \cong \pi_2(\text{Diff}(\BT))$ which is rationally trivial by Theorem \ref{htp-g-Top-Diff}. Hence we can apply Proposition \ref{modify-nh}, up to change $f_i$ by $Nf_i$ for some suitable $N$, we can assume $F\mid_{\partial D^4 =S^3}$ represents a trivial element in $\pi_3(\text{BDiff}(\BT))$. Thus we can extend $F$ to $D^4$ and get a map $F$ from $M$ to $\text{BDiff}(\BT) $.

We are left to show the total space of the bundle  we have just constructed has finite abelian fundamental group. Note that for this, we only need the information of the map $\tilde{F} : M \rightarrow \text{BG}(\BT)$. In fact, since we can assume that the diffeomorphism we used here lies in the identity component of $\text{Diff}(\BT)$, we can assume $\tilde{F} : M \rightarrow \text{BG}_0(\BT)$, where $\text{G}_0(\BT)$ is the identity component of $\text{G}(\BT)$. Now by Lemma \ref{htp-g-g}, we have $\text{BG}_0(\BT)$ is homotopy equivalent to $K(\BZ^{n},2)$, hence the maps from $M$ to $\text{BG}_0(\BT)$ are classified by $H^2(M, \BZ^{n})$. Since our $M$ is a simply connected $4$-manifold, $H^2(M, \BZ^{n}) =  Hom(H_2(M),\BZ^{n})$. And we know precisely what this map is by reading off information from $F\mid_{\bigvee_{i=1}^{n} S^2}$. Denote the corresponding map in $Hom(H_2(M),\BZ^{n})$ by $\alpha$, we have the image of $\alpha$ in $\BZ^{n}$ has finite index.

On the other hand, we have a universal fibration $\eta$ over $\text{BG}_0(\BT)$ with fiber the $n$-dimensional torus $T$ and the bundle over $M$ we constructed as a fibration is the same as pullback the fibration $\eta$ via $\tilde{F}$. Denote the total space of $\eta$ by $E_\eta$.

\textbf{Claim.} The total
space  $E_\eta$ is simply connected.

Assuming the claim, we have the following two  long exact sequences of homotopy groups

\begin{tikzcd}
\dots \arrow[r]  & \pi_2 (M)\arrow[r,"\partial"] \arrow[d, "\alpha_{2\ast}"] & \pi_1(\BT) \arrow[r]\arrow[d,"\alpha_{1\ast}"] & \pi_1(E) \arrow[d] \arrow[r] &  \pi_1 (M)\arrow[r] \arrow[r] \arrow[d] & \dots\\
\dots \arrow[r]  & \pi_2 (\text{BG}_0(\BT)) \arrow[r,"\partial_{\eta}"] & \pi_1(T) \arrow[r] & \pi_1(E_\eta) \arrow[r] &  \pi_1 (\text{BG}_0(\BT))\arrow[r] &\dots
\end{tikzcd}

Note that the map $\partial_{\eta}$ on the lower sequence is surjective  by the claim. On the other hand $\pi_2 (\text{BG}_0(\BT)) \cong \pi_1(T) \cong \BZ^n$, we have $\partial_{\eta}$ is an isomorphism. The induced $\alpha_{1\ast}: \pi_1(\BT) \rightarrow \pi_1(T)$ is also an  isomorphism. Moreover, we have the map $\alpha_{2\ast} : \pi_2 (M)\rightarrow \pi_2 (\text{BG}_0(\BT))$ is determined by the map $\alpha$ and the image of it in $\text{BG}_0(\BT)$ has finite index. Hence,  we have the image of $\partial$ in $\pi_1(\BT)$ is a finite index subgroup. Since $M$ is simply connected, we have $\pi_1(E)$ is a finite abelian group. 

\textbf{Proof of  the Claim.} Note first that the total space of the universal principal bundle for the torus  $T$ (considered as a Lie group) is contractible with $BT = K(\BZ^n,2)$. As a $T$-fibration, it is classified by a map $\beta:\text{B}T \rightarrow \text{BG}_0(T) = \text{BG}_0(\BT)$. Just as before, we have the following two  long exact sequences of homotopy groups.

\begin{tikzcd}
\dots \arrow[r]  & \pi_2 (\text{B}T) \arrow[r,"\partial_\beta"] \arrow[d, "\beta_{2\ast}"] & \pi_1(T) \arrow[r]\arrow[d,"\beta_{1\ast}"] & \pi_1(\ast) \arrow[d] \arrow[r] &  \pi_1 (\text{B}T)\arrow[r] \arrow[r] \arrow[d] & \dots\\
\dots \arrow[r]  & \pi_2 (\text{BG}_0(\BT)) \arrow[r,"\partial_{\eta}"] & \pi_1(T) \arrow[r] & \pi_1(E_\eta) \arrow[r] &  \pi_1 (\text{BG}_0(\BT))\arrow[r] &\dots
\end{tikzcd}

Now the map $\partial_\beta$ and $\beta_{1\ast}$ are  isomorphisms, we have the map $\partial_\eta$ is a surjection. In fact since $\pi_2 (\text{BG}_0(\BT)) \cong \pi_1(T) \cong \BZ^n$, we have $\partial_\eta$ must be an isomorphism. Now the claim follows from the fact that $\pi_1(\text{BG}_0(\BT)) =0$.

\bibliographystyle{amsplain}

\end{document}